\documentclass[preprint,sort&compress,12pt]{elsarticle}
\usepackage{amsmath}
\usepackage{mathrsfs}
\usepackage{stmaryrd}
\usepackage{mathrsfs}
\usepackage{bm}
\usepackage{amsmath}
\usepackage{amsfonts}
\usepackage{amssymb}
\usepackage{amsthm}

\journal{\quad }
\newcommand{\mf}{\mathbf}
\newcommand{\mm}{\mathrm}
\begin{document}
\begin{frontmatter}
\title{A Note on Large Time Behavior of Velocity in the Baratropic Compressible Navier--Stokes Equations \tnoteref{S}}
\author[FZ,FJ]{Fei Jiang\corref{cor1}}
\ead{jiangfei0591@163.com} \cortext[cor1]{Corresponding author: Tel
+86 15001201710.}
\address[FZ]{College of Mathematics and Computer Science, Fuzhou University, Fuzhou, 361000, China.}
\address[FJ]{Institute of Applied Physics and Computational Mathematics, Beijing, 100088, China.}
\begin{abstract}
Recently, for periodic initial data with initial density allowed to
vanish, Huang and Li \cite{HXDLJEh2} establish the global existence
of strong and weak solutions for the two-dimensional compressible
Navier--Stokes equations with no restrictions on the size of initial
data provided the bulk viscosity coefficient is $\lambda=\rho^\beta$
with $\beta> 4/3$. Moreover, the large-time behavior of the strong
and weak solutions are also obtained, in which the  velocity
gradient strongly converges to zero in $L^2$ norm. In this note, we
further point out that the velocity strongly converges to an
equilibrium velocity in $H^1$ norm, in which the equilibrium
velocity is uniquely determined by the initial data. Our result can
also be regarded a correction for the result of large-time behavior
of velocity in \cite{PMOS311}. \\[-1.5em]
\end{abstract}

\end{frontmatter}


\newtheorem{thm}{Theorem}[section]
\newtheorem{lem}{Lemma}[section]
\newtheorem{pro}{Proposition}[section]
\newtheorem{cor}{Corollary}[section]
\newproof{pf}{Proof}
\newdefinition{rem}{Remark}[section]
\newtheorem{definition}{Definition}[section]
\section{Introduction}
\label{Intro} \numberwithin{equation}{section}

In this note, we are concerned with the two-dimensional barotropic
compressible Navier--Stokes equations which read as follows:
\begin{eqnarray}
&&   \label{0101}                  \partial_t\rho+ \mathrm{div}(\rho{\mathbf{v}})=0,  \\
  &&   \label{0102}         \partial_t(\rho {\mathbf{v}})+\mathrm{div}(\rho{\mathbf{v}}\otimes {\mathbf{v}})+\nabla P(\rho)=
                    \mu \Delta{\mathbf{v}}+\nabla ((\lambda+\mu)
                    \mathrm{div}{\mathbf{v}}),
 \end{eqnarray}where $\rho$ and $\mf{v}$
represent the density and velocity respectively, and the pressure
$P$ is given by\begin{equation*} P(\rho)=a\rho^\gamma,\quad
\gamma>1.\end{equation*}Here $a=e^S>0$ is the constant determined by
the entropy constant $S$, and $\gamma\geq 1$
 the adiabatic constant. Values of $\gamma$ have their
own physical significance, and are also take important part in the
existence of solutions (see \cite{CS,WWWJFGZSM3,JFTZON6310,GRJFYJ}
for example). The viscosity coefficients satisfy the following
hypothesis:
\begin{equation*}\mu=\mbox{constant},\ \lambda(\rho)=b\rho^\beta,\  b>0,\ \beta> 0.
\end{equation*} As in \cite{HXDLJEh2}, we consider the Cauchy problem
with the given initial density  $\rho_0$ and the given initial
momentum $\mf{m}_0$, which are periodic with period 1 in each space
direction $x_i$, $i=1$, $2$, i.e., functions defined on
$\mathbb{T}^2=\mathbb{R}^2/\mathbb{Z}^2$. We require that
\begin{equation*}\rho(\mathbf{x},0)=\rho_0(\mathbf{x}),\
\rho\mathbf{v}(\mf{x},0)=\mf{m}_0(\mf{x}),\ \mf{x}=(x_1,x_2)\in
\mathbb{T}^2.
\end{equation*}

There is a huge literature concerning the theory of strong and weak
solutions for the system of the multidimensional compressible
Navier--Stokes with constant viscosity coefficients. The local
existence and uniqueness of classical solutions are known in
\cite{NJLB91,SJOSA3} in the absence of vacuum and recently, for
strong solutions also, in \cite{UCYCHJKHJ8,SCHJKHJ15,SMTRSMQR} for
the case that the initial density need not be positive and may
vanish in open sets. The global classical solutions were first
obtained by Matsumura and Nishida \cite{MANTTJ601} for initial data
close to a non-vaccum equilibrium in some Sobolev space $H^s$.
Later, Hoff \cite{GHDJ12} studied the problem for discontinuous
initial data. For the existence of solutions for large data, the
major breakthrough is due to Lions \cite{LPLMTFM98} (see also
Feireisl \cite{FETD2,FENAPHOFJ35801}), where he obtained global
existence of weak solutions, defined as solutions with finite
energy, when the exponent $\gamma$ is suitably large. The main
restriction on initial data is that the initial energy is finite, so
that the density is allowed to vanish initially. Recently, Huang, Li
and Xin \cite{HXLJXZPG} established the global existence and
uniqueness of classical solutions to the Cauchy problem for the
isentropic compressible Navier--Stokes equations in the
three-dimensional space with smooth initial data which are of small
energy but possibly large oscillations; in particular, the initial
density is allowed to vanish, even has compact support.

However, there are few results regarding global strong solvability
for equations of multi-dimensional motions of viscous gas with no
restrictions on the size of initial data. One of the first ever ones
is due to Vaigant--Kazhikhov \cite{VVAKAVOS3} who obtained a
remarkable result which can be stated that the two-dimensional
system (\ref{0101})--(\ref{0102}) admits a unique  global strong
solution for large initial data away from vacuum provided $\beta>3$.
Lately, Perepelitsa \cite{PMOS311} proved the global existence of a
weak solution with uniform lower and upper bounds on the density, as
well as the decay of the solution to an equilibrium state in a
special case that
\begin{equation*}\beta>3,\ \gamma=\beta.
\end{equation*}
when the initial density is away from vacuum. Very recently, Jiu,
Wang and Xin \cite{JQWYXZPG} consider classical solutions and
removed the condition that the initial density should be away from
vacuum in Vaigant-Kazhikhov \cite{VVAKAVOS3} but still under the
same condition that $\beta>3$ as that in \cite{VVAKAVOS3}. No long
after, Huang and Li establish the global existence of strong and
weak solutions provided $\beta> 4/3$ and $\gamma>1$.

 Before stating the exciting
result of Huang and Li, we explain the notations and conventions
used throughout this paper. We denote
\begin{equation*}\int f\mm{d}\mf{x}=\int_{\mathbb{T}^2}f\mm{d}\mf{x},\ \bar{f}=\frac{1}{|\mathbb{T}^2|}\int f\mm{d}\mf{x}.
\end{equation*}
For $1\leq r\leq \infty$, we also denote the standard Lebesgue and
Sobolev spaces as follows:
\begin{equation*}L^r=L^r(\mathbb{T}^2),\ W^{s,r}=W^{s,r}(\mathbb{T}^2),\ H^s=W^{s,2}.
\end{equation*}
Then, we state the Huang and Li's result concerning the global
existence and large-time behavior of strong solutions as follows:
\begin{thm}\label{thm:0101} Assume that
\begin{equation}\label{0108}\beta>4/3,\ \gamma>1,
\end{equation}and that the initial data $(\rho_0,\mf{m}_0)$ satisfy
that for some $q>2$,
\begin{equation*}\label{0109}0\leq \rho_0\in W^{1,q},\ \bar{\rho}_0> 0,\ \mf{v}_0\in H^1,\
\mf{m}_0=\rho_0\mf{v}_0.
\end{equation*}
Then the problem (\ref{0101})--(\ref{0102}) has a unique global
strong solution $(\rho,\mf{v})$ satisfying
\begin{equation*}\label{0110}\left\{\begin{array}{l}
  \rho\in C([0,T],W^{1,1}),\ \rho_t\in L^\infty(0,T; L^2), \\
  \mf{v}\in L^\infty(0,T;H^1)\cap L^{(q+1)/q}(0,T;W^{2,q}), \\
  t^{1/2}\mf{v}\in L^2(0,T;W^{2,q}),\ t^{1/2}\mf{v}_t\in L^2(0,T;H^1), \\
  \rho\mf{v}\in C([0,T],L^2), \ \sqrt{\rho}\mf{v}_t\in
  L^2(\mathbb{T}^2\times(0,T)),
\end{array}\right.\end{equation*}for any $0<T<\infty$. Moreover, if
\begin{equation}\label{0111}\beta>3/2,\ 1<\gamma<3(\beta-1),\end{equation}
there exists a constant $C$ independent of $T$ such that
\begin{equation}\label{0112}\sup_{0\leq t\leq T}\|\rho(\cdot,t)\|_{L^\infty}\leq C,\end{equation}
\begin{equation}\label{0113}\sup_{0\leq t\leq T}\|\mf{v}(\cdot,t)\|_{H^1}\leq C,\end{equation}
and the following large-time behavior holds:
\begin{equation}\label{0114}\lim_{t\rightarrow \infty}(\|\rho-\bar{\rho}_0\|_{L^p}+\|\nabla \mathbf{v}\|_{L^2})=0,\end{equation}
for any $p\in [1,\infty)$.
\end{thm}
\begin{rem}\label{rem:0101} The results above can be found in \cite[Theorem 1.1]{HXDLJEh2}, except for the  estimate (\ref{0113}).
Fortunately we can obtain (\ref{0113}) by \cite[Proposition
3.5]{HXDLJEh2}.
\end{rem}

The result (\ref{0114}) above indicates that the density $\rho(t)$
strongly converges to the equilibrium density $\bar{\rho}_0$ in
$L^p$ norm as $t\rightarrow \infty$. We naturally propose an
interesting  question of whether there exists an equilibrium
velocity $\mf{v}_s$ such that the velocity $\mathbf{v}(t)$ strongly
converges $\mf{v}_s$ in some norm  as $t\rightarrow \infty$. In this
note, we give the positive result. Next, we state our result, which
will be proved in Section 2.
\begin{thm}\label{thm:0102} Assume that the strong solution $(\rho,\mf{v})$
is provided by Theorem \ref{thm:0101}. If  $(\rho,\mf{v})$ satisfies
(\ref{0113}) and (\ref{0114}), then
\begin{equation}\label{0115}\lim_{t\rightarrow \infty}\| \mathbf{v}-\mathbf{v}_s\|_{H^1}=0,\end{equation}
where
\begin{equation}\label{0116}{{\mathbf{v}}}_s:=\frac{1}{\bar{\rho}_0|{\mathbb{T}^2}|}\int
\rho_0{\mathbf{v}}_0\mm{d}\mf{x}\mbox{ is a constant vector.
}\end{equation}
\end{thm}

\begin{rem} Assume that (\ref{0108}) holds and that the initial data
$(\rho_0, \mf{m}_0)$ satisfies that $0\leq \rho_0\in L^\infty$,
$\mf{v}_0\in H^1$, $\mf{m}_0=\rho_0\mf{v}_0$. Then the problem
(\ref{0101})--(\ref{0102}) possesses
 at least one global weak solution $(\rho,\mf{v})$. Moreover, if
$\beta$ and $\gamma$ satisfy (\ref{0111}), there exists a constant
$C$ independent of $T$ such that (\ref{0112})--(\ref{0114}) hold
true (see \cite[Theorem 1.2]{HXDLJEh2}). We mention that such weak
solution also satisfies (\ref{0115}).
\end{rem}
\begin{rem} When the initial density is away from vacuum, Perepelitsa \cite{PMOS311} proved the global existence of a
weak solution, as well as the convergence of the solution to an
equilibrium state in a special case that $\beta>3$ and
$\gamma=\beta$, where the author considered that equilibrium
velocity is zero vector. According to Theorem \ref{thm:0102},
 the equilibrium velocity is uniquely determined by the relation
(\ref{0116}), and is not zero vector in the general case.
\end{rem}
\section{Proof of Theorem \ref{thm:0102}}
In this section, we start to prove   Theorem \ref{thm:0102}. First,
exploiting (\ref{0113}) and the fact that $H^1\hookrightarrow L^2$
is compact, we have  that for any sequence
$\{t_n\}_{n=1}^\infty\subset (0,\infty)$, there exists a subsequence
$\{t_{n_m}\}_{m=1}^\infty\subset \{t_n\}_{n=1}^\infty$, such that
\begin{eqnarray}\label{0201}
&&\mf{v}(t_{n_m})~{\rightharpoonup}~{\mf{v}}_M\mbox{ weakly in
}H^1,\\
&&\label{n0201}\mf{v}(t_{n_{m}})~{\rightarrow}~{\mf{v}}_M\mbox{
strongly in
}L^2,\\
&& n_m\rightarrow \infty\mbox{ as }m\rightarrow \infty.\nonumber
 \end{eqnarray}
Thanks to the condition (\ref{0114}), we see that
$\lim_{t\rightarrow \infty}\|\nabla \mf{v}(t)\|=0$, so ${\mf{v}}_M$
must be a constant
 vector.

Next we shall show that the constant vector ${\mf{v}}_M$ doses not depend
on the particular choice of subsequences. To this end, integrating
the equation (\ref{0102}), we can deduce the momentum conservation
\begin{eqnarray}\label{0202}
\int\rho(t)\mf{v}(t)\mm{d}\mf{x}=\int\rho_0\mf{v}_0\mm{d}\mf{x}.
 \end{eqnarray}
Letting  $t:=t_m\rightarrow \infty$ in (\ref{0202}), and using
(\ref{0114}) and (\ref{n0201}), then we obtain
\begin{eqnarray*}
\int\bar{\rho}_0{\mf{v}}_M\mm{d}\mf{x}=\lim_{t_m\rightarrow
\infty}\int\rho(t_m)\mf{v}(t_m)\mm{d}\mf{x}=\int\rho_0\mf{v}_0\mm{d}\mf{x},
 \end{eqnarray*}
which  yields
\begin{eqnarray}
{\mf{v}}_{M}\equiv {\mf{v}}_s:=
\frac{\int\rho_0\mf{v}_0\mm{d}\mf{x}}{\bar{\rho}_0|\mathbb{T}^2|}.
 \end{eqnarray}

Consequently, making use of the convergence of velocity in
(\ref{0114}) and (\ref{n0201}), we can conclude that for any
sequence $\{t_n\}_{n=1}^\infty\subset (0,\infty)$, there exists a
subsequence $\{t_{n_m}\}_{m=1}^\infty\subset \{t_n\}_{n=1}^\infty$,
such that
\begin{eqnarray}\label{0205}
&& \mf{v}(t_{n_{m}})~{\rightarrow}~{\mf{v}}_s\mbox{ strongly in }H^1
\mbox{ as }m\rightarrow \infty.
 \end{eqnarray}
Hence (\ref{0115}) holds, since the sequence
$\{t_n\}_{n=1}^\infty\subset (0,\infty)$ is arbitrary. This
completes the proof of Theorem \ref{thm:0102}.

\newcommand\ack{\section*{Acknowledgement}}
\newcommand\acks{\section*{Acknowledgements}}
\renewcommand\refname{References}
\renewenvironment{thebibliography}[1]{%
\section*{\refname}
\list{{\arabic{enumi}}}{\def\makelabel##1{\hss{##1}}\topsep=0mm
\parsep=0mm
\partopsep=0mm\itemsep=0mm
\labelsep=1ex\itemindent=0mm
\settowidth\labelwidth{\small[#1]}%
\leftmargin\labelwidth \advance\leftmargin\labelsep
\advance\leftmargin -\itemindent
\usecounter{enumi}}\small
\def\newblock{\ }
\sloppy\clubpenalty4000\widowpenalty4000
\sfcode`\.=1000\relax}{\endlist}
%


\bibliographystyle{model1-num-names}







\end{document}